\documentclass[11pt]{article}
\usepackage[letterpaper,hmargin=1in,vmargin=1.25in]{geometry}
\usepackage{latexsym, amsfonts, amssymb, amsmath, amsthm, graphicx}

\newcommand{\charf}[1]{\mbox{\raise.48ex\hbox{$\chi$}$_{#1}$}}
\def\~{\tilde }
\def\C{{\cal C}}

\def\R{{\mathbb R}}

\def\Q{{\sf Q}}
\def\P{{\sf P}}
\def\I{{\rm I}}
\def\E{{\sf E}}

\def\1{{\sf 1}}
\def\K{{\sf K}}

\def\var{{\rm Var}}

\theoremstyle{plain}
\newtheorem{theorem}{Theorem}[section]
\newtheorem{lemma}[theorem]{Lemma}
\newtheorem{corollary}[theorem]{Corollary}

\theoremstyle{definition}

\newtheorem{example}[theorem]{Example}

\begin{document}

\title{The simple random walk and max-degree walk on a directed graph}

\author {Ravi Montenegro \thanks{Department of Mathematical Sciences, University of
      Massachusetts Lowell, Lowell, MA 01854, ravi\_montenegro@uml.edu.}
}

\date{}

\maketitle
 
\begin{abstract}
\noindent
We show bounds on total variation and $L^{\infty}$ mixing times, spectral gap and magnitudes of the complex valued eigenvalues of a general (non-reversible non-lazy) Markov chain with a minor expansion property. This leads to the first known bounds for the non-lazy simple and max-degree walks on a (directed) graph, and even in the lazy case they are the first bounds of the optimal order. In particular, it is found that within a factor of two or four, the worst case of each of these mixing time and eigenvalue quantities is a walk on a cycle with clockwise drift.

\vspace{1ex}
\end{abstract}

\noindent
\section{Introduction}

Markov chains are a key tool in approximation algorithms for combinatorial counting problems and for sampling from discrete spaces. Surprisingly, little is known about the convergence rate of a Markov chain with no holding probability. Even for the lazy simple random walk (i.e. strongly aperiodic nearest neighbor walk) on an undirected graph the order of magnitude for the slowest mixing time is not known. 

More specifically, consider an undirected graph with $m$ edges, $n$ vertices and maximum degree $d$. The lazy simple random walk is known to converge in $O(m^2\log(m/\epsilon))$ steps, and so the lazy max-degree walk mixes in time $O(n^2d^2\log(nd/\epsilon))$ as well. However, there are no examples for which either bound is of the correct order. We solve this riddle by giving new bounds for simple and max-degree walks which are better then these. Moreover, our bounds apply to directed graphs, require no holding probability, and are nearly sharp.

To state our results, recall that an Eulerian graph is a strongly connected directed graph such that each vertex has the same in and out-degrees. This is the natural directed analog of an undirected graph, as any undirected graph can be made into an Eulerian graph by replacing each undirected edge with two directed edges. Two natural walks on a graph will be considered. For the simple random walk choose a neighbor uniformly at random and go there, while in the max-degree walk choose a neighbor with probability $1/d$ each and otherwise do nothing. 

Which (non-lazy) directed walks mix rapidly? Certainly the walk should not get stuck drifting between sets of equal sizes, such as from one bipartition to another (e.g. simple walk on a cycle with an even number of vertices). To avoid this it suffices that if the walk starts in a set of size $x\leq 1/2$, then the set of adjacent vertices has size $>x$. For instance, a max-degree walk on a strongly connected graph with a self-loop at each vertex. 

We now give our main results. Note that $\tau(\epsilon)$ is total variation mixing time (time to converge at an ``average'' vertex), $\tau_{\infty}(\epsilon)$ is $L^{\infty}$ mixing time (time to converge at every vertex), $\lambda_i\neq 1$ is any non-trivial (complex-valued) eigenvalue of the transition matrix, $\lambda$ is the spectral gap, and $N(A)=\{y\in V:\,\exists x\in A,\,\P(x,y)>0\}$ is the neighborhood of set $A$.

\begin{corollary} \label{cor:Eulerian-simple}
The simple random walk on an Eulerian graph with $m$ edges satisfies
$$
\lambda 
   \geq 1-\cos \frac{2\pi}{m}
   \approx \frac{2\pi^2}{m^2}\,.
$$
If it satisfies the expansion condition that
$$
\forall A\subset V,\,\pi(A)\leq 1/2,\,\forall v\in V:\,\pi\left(N(A)\setminus v\right) \geq \pi(A)
$$
then also
\begin{eqnarray*}
1-|\lambda_i|
   &\geq& 1-\cos\frac{2\pi}{m}
   \approx \frac{2\pi^2}{m^2} \\
\tau(\epsilon) 
   &\leq& \frac{1}{-\log \cos\frac{2\pi}{m}}\log \frac{1-2/m}{\epsilon} 
   \approx \frac{m^2}{2\pi^2}\log\frac{1}{\epsilon} \\
\tau_{\infty}(\epsilon)
   &\leq& \min\left\{\frac{\log \frac{m-2}{2}+\log\frac{1}{\epsilon}}{-\log \cos\frac{2\pi}{m}},\,\frac{m^2}{6}+\frac{m^2}{8}\log\frac{1}{\epsilon}\right\}\,\delta_{\epsilon\leq 1}+m^2\frac{1+3\epsilon}{3(1+\epsilon)^3}\,\delta_{\epsilon>1} \\
   &\approx& \frac{m^2}{2\pi^2}\log \frac{m-2}{2\epsilon}\,\delta_{\epsilon< 1/m} + \frac{m^2}{8}\log\frac{4}{\epsilon}\,\delta_{\epsilon\in[1/m,1]} + \frac{m^2}{(1+\epsilon)^2}\,\delta_{\epsilon>1}
\end{eqnarray*}

For the lazy simple random walk the bound on $\lambda$ is a factor two smaller, the expansion condition is replaced by strong connectivity, and in the remaining bounds replace $m$ by $2m$.
\end{corollary}

It follows that every lazy simple Eulerian walk converges in the same $\tau(\epsilon)=O(m^2\log(1/\epsilon))$ steps required for a cycle walk, improving on and generalizing the classical result $\tau(\epsilon)=O(m^2\log(m/\epsilon))$ for a lazy simple undirected walk. This can be further improved on by an order of magnitude in the special case of a walk on a regular graph, or equivalently of a max-degree walk.

\begin{corollary} \label{cor:Eulerian-max}
The max-degree walk on an Eulerian graph with $n$ vertices and max-degree $d$ satisfies
$$
\lambda
   \geq \frac 2d\left(1-\cos \frac{\pi}{n}\right)
   \approx \frac{\pi^2}{n^2 d}\,.
$$
If it satisfies the expansion condition that
$$
\forall A\subset V,\,|A| \leq |V|/2:\,|N(A)| > |A|
$$
then also
\begin{eqnarray*}
1-|\lambda_i| 
   &\geq& \frac 2d\left(1-\cos \frac{\pi}{n}\right)
   \approx \frac{\pi^2}{n^2 d} \\
\tau(\epsilon) 
   &\leq& \frac{1}{-\log (1-\frac 2d(1-\cos \frac{\pi}{n}))}\log \frac{1-1/n}{\epsilon} 
   \approx \frac{n^2 d}{\pi^2}\log\frac{1}{\epsilon} \\
\tau_{\infty}(\epsilon)
   &\leq& \min\left\{\frac{\log(n-1)+\log\frac{1}{\epsilon}}{-\log (1-\frac 2d(1-\cos \frac{\pi}{n}))},\,\frac{n^2d}{3}+\frac{n^2d}{4}\log\frac{1}{\epsilon}\right\}\,\delta_{\epsilon\leq 1}+n^2d\,\frac 23\,\frac{1+3\epsilon}{(1+\epsilon)^3}\,\delta_{\epsilon>1} \\
   &\approx& \frac{n^2 d}{\pi^2}\log\frac{n-1}{\epsilon}\,\delta_{\epsilon<1/n} + \frac{n^2 d}{4}\log\frac{4}{\epsilon}\,\delta_{\epsilon\in[1/n,1]} + \frac{2\,n^2 d}{(1+\epsilon)^2}\,\delta_{\epsilon>1}
\end{eqnarray*}

For the lazy max-degree walk the bounds on $\lambda$ is a factor two smaller, the expansion condition is replaced by strong connectivity, and in the remaining bounds replace $d$ by $2d$.
\end{corollary}

Just how good are these bounds? For the simple random walk on the cycle with an odd number of vertices $n$ (so $m=2n$ and $d=2$) the spectral gap bound is off by a factor of $4$, the eigenvalue bounds are exact, and the upper bounds on $\tau(\epsilon)$ become lower bounds if $\log\frac{1-\pi_*}{\epsilon}$ is replaced by $\log\frac{1}{2\epsilon}$ (where $\pi_*=2/m$ and $\pi_*=1/n$ respectively). More generally, we define a precise notion of rate of expansion, and show that a cycle walk with clockwise drift will be within a factor two of being the slowest mixing, not only among simple or max-degree walks, but among all Markov chains with this rate of expansion!

An interesting aspect of our argument is that it uses the Evolving set methodology of Morris and Peres \cite{MP05.1}, in an improved form given by this author \cite{Mon06.1} which bounds total variation distance directly, without going through $L^2$ distance. Related bounds also show that with relative entropy and $L^2$ mixing times the cycle walk is again asymptotically nearly the slowest walk.

The paper proceeds as follows. In Section \ref{sec:preliminaries} we review the Evolving set methodology. This is followed in Section \ref{sec:general} with a proof of our main mixing result, a generalization of the simple and max-degree Eulerian walks considered above. In Section \ref{sec:distances} this is extended to a bound on convergence rates in distances other than total variation. The Appendix contains proofs of inequalities used in showing our results.


\section{Review of Mixing and Evolving Sets} \label{sec:preliminaries}

We begin by reviewing mixing time theory, and particularly Evolving Set ideas.

Let $\P$ be a finite irreducible Markov kernel on state space $V$ with stationary distribution $\pi$, that is, $\P$ is a $|V|\times|V|$ matrix with entries in $[0,1]$, row sums are one, $V$ is connected under $\P$ ($\forall x,y\in V\,\exists t:\,\P^t(x,y)>0$), and $\pi$ is a distribution on $V$ with $\pi\P=\pi$. The time-reversal $\P^*$ is given by $\P^*(x,y)=\frac{\pi(y)\P(y,x)}{\pi(x)}$ and has stationary distribution $\pi$ as well. If $A,B\subset V$ the ergodic flow from $A$ to $B$ is given by $\Q(A,B)=\sum_{x\in A,y\in B} \pi(x)\P(x,y)$. Given initial distribution $\sigma$, the $t$-step discrete time distribution is given by $\sigma\P^t$. 

If the walk is strongly connected and aperiodic then $\sigma\P^t\xrightarrow{t\to\infty}\pi$. Our goal is to measure the rate of convergence. One of the more widely used notions of distance between distributions is the variation distance, 
$$
\|\sigma-\pi\|_{TV} = \frac 12\sum_{x\in V} |\sigma(x)-\pi(x)|\,.
$$
The mixing time $\tau(\epsilon)$ denotes the worst-case number of steps required for the total variation distance $\|\P^t(x,\cdot)-\pi\|_{TV}$ to drop to $\epsilon$.

Many bounds on mixing time are shown by working with the spectral gap, which is just the gap between the two largest eigenvalues of the walk $\frac{\P+\P^*}{2}$, that is,
$$
\lambda = \min_{i\neq 0} 1-\lambda_i\left(\frac{\P+\P^*}{2}\right)
 = \inf_{\var(f)\neq 0} \frac{\frac 12\sum_{x,y\in V}(f(x)-f(y))^2\pi(x)\P(x,y)}{\frac 12\sum_{x,y\in V} (f(x)-f(y))^2\pi(x)\pi(y)}
$$
where $\{\lambda_i(\K)\}$ denotes the eigenvalues of Markov chain $\K$, and $\lambda_0(\K)=1$. 

Our results are based on a theorem of Montenegro \cite{Mon06.4}, which is in turn proven by working with a lower bound on mixing given in \cite{Mon06.2}, and with the Evolving set methodology of Morris and Peres \cite{MP05.1} in a stronger form given by Montenegro \cite{Mon06.1}.

\begin{theorem} \label{thm:main}
Consider a finite, irreducible Markov chain. Given $f:[0,1]\to\R_+$ such that $\forall a\in(0,1/2]:\,0<f(a)\leq f(1-a)$, let the $f$-congestion $\C_f=\max_{A\subset V,\,\pi(A)\leq 1/2}\C_f(A)$, where
$$
\C_f(A) = \frac{\int_0^1 f(\pi(A_u))\,du}{f(\pi(A))}
\quad and\quad 
A_u = \{v \in V:\,\Q(A,v)\geq u\pi(v)\}\,.
$$
Then, the $t$-step Markov chain satisfies
$$
\frac 12\,|\lambda_i|^t \leq \max_{x\in V} \|\P^t(x,\cdot)-\pi\|_{TV} \leq \left(\max_{\pi(A)\leq 1/2} \frac{\pi(A)\pi(A^c)}{f(\pi(A))}\right) \,\left(\max_{x\in V}\frac{f(\pi(x))}{\pi(x)}\right)\,\C_f^t\,,
$$
while every (complex valued) eigenvalue $\lambda_i\neq 1$ of the transition probability matrix satisfies 
$$
1-|\lambda_i| \geq 1-\C_f\,.
$$

If instead $f$ satisfies the weaker condition $\forall a\in(0,1):\,f(a)>0$, then the result still holds, but with $\C_f=\max_{A\subset V}\C_f(A)$.
\end{theorem}

The proof is not yet in print, and so it is included in the Appendix to this paper. 

The $f$-congestion $\C_f(A)$ is a measure of the expansion (or congestion) of a random walk, and seems particularly well suited towards showing geometric bounds on mixing times. While it is generally not easy to calculate directly, the following lemma of \cite{Mon06.1} makes it possible to bound $\C_f$ in terms of isoperimetric quantities.

\begin{lemma} \label{lem:worst-case}
Given a concave function $f:[0,1]\rightarrow \R$ and
two non-increasing functions $g,\,\hat{g}:\,[0,1]\rightarrow[0,1]$
such that $\int_0^1 g(u)\,du=\int_0^1 \hat{g}(u)\,du$ and
$\forall t\in[0,1]:\,\int_0^t g(u)\,du \geq \int_0^t \hat{g}(u)\,du$, 
then
$$
\int_0^1 f\circ g(u)\,du \leq \int_0^1 f\circ \hat{g}(u)\,du\,.
$$
\end{lemma}

Our interest is in bounding $\int_0^1 f(\pi(A_u))\,du$ for a concave function $f$. To apply the lemma we minimize the integral $\int_0^t \pi(A_u)\,du$, recalling that $\pi(A_u)$ is a decreasing function of $u$ and noting that $\int_0^1 \pi(A_u)\,du=\pi(A)$, while taking into account whatever constraints are given by the problem of interest. 
A particularly useful constraint to consider is the {\em modified ergodic flow}
$$
\forall A\subset V:\,\Psi(A) = \frac 12\int_0^1|\pi(A_u)-\pi(A)|\,du\,.
$$
By the property $\int_0^1\pi(A_u)\,du=\pi(A)$ it follows that $\Psi(A)$ is the area below $\pi(A_u)$ and above $\pi(A)$, while also the area below $\pi(A)$ and above $\pi(A_u)$. See Figures \ref{fig:uniform} and \ref{fig:general} for two such examples.

A useful interpretation of $\Psi(A)$ is as the smallest ergodic flow from set $A$ to a set $B$ of size $\pi(B)=\pi(A^c)$ (see \cite{Mon06.1}), that is
$$
\Psi(A) = \min_{\substack{B\subset V,\,v\in V,\\ \pi(B)\leq\pi(A^c)<\pi(B\cup v)}} \Q(A,B)+\frac{\pi(A^c)-\pi(B)}{\pi(v)}\,\Q(A,v)\,.
$$
When the distribution $\pi$ is uniform then this simplifies to
$$
\Psi(A) = \min_{\substack{B\subset V,\\ \pi(B)=\pi(A^c)}} \Q(A,B)\,,
$$
while in general if the walk is lazy (i.e. $\forall x\in V:\,\P(x,x)\geq 1/2$) then $\Psi(A)=\Q(A,A^c)$ with the worst set $B$ being $B=A^c$.

Now, what is a good choice of function $f$ for the $f$-congestion? In \cite{Mon06.4} it was suggested that if the modified ergodic flow $\Psi(A)\geq C,\,\forall A\subset V$ for some constant $C$ not depending on set size, such as the simple random walk on an odd length cycle (with $\Psi(A)\geq 1/m$), then it is best to work with $f(a)=\sin(\pi\,a)$. In this case Theorem \ref{thm:main} implies
$$
\|\P^t(x,\cdot)-\pi\|_{TV} \leq (1-\pi_*)\C_{\sin(\pi\,a)}^t
\quad{\rm and}\quad
1-|\lambda_i| \geq 1-\C_{\sin(\pi\,a)}\,,
$$
where $\pi_*=\min_{v\in V}\pi(v)$.


\section{General random walks} \label{sec:general}

We now set out to show our main result, eigenvalue and total variation mixing bounds for general random walks (the $L^{\infty}$ case will be dealt with in the next section). Two corollaries of this will be the specific walks on Eulerian graphs discussed in the introduction. In particular, we will find that even when general Markov chains are considered, a walk with clockwise drift on a cycle is still within a factor two of being the slowest mixing Markov chain.

Recall from the preliminaries that $\Psi(A)$ is both a measure of area between the curves $y=\pi(A_u)$ and $y=\pi(A)$, and also a measure of the worst ergodic flow from $A$ into a set of size $\pi(A^c)$. While this will play a key role in our proof, our main theorem will involve a slightly weaker quantity. In practice these two will usually be the same. Given $A\subset V$, let $\hat\Q(A,x)=\min\{\Q(A,x),\pi(x)/2\}$ and define
\begin{eqnarray*}
\hat\Psi(A)
  &=& \min_{\substack{B\subset V,\,v\in V,\\ 
             \pi(B)\leq\pi(A^c)<\pi(B\cup v)}} 
   \sum_{x\in B} \hat\Q(A,x) + \frac{\pi(A^c)-\pi(B)}{\pi(v)}\hat\Q(A,v)\,.
\end{eqnarray*}
As with $\Psi(A)$, when the distribution is uniform then this is just $\hat\Psi_{min}=\min_{\pi(B)=\pi(A^c)}\hat\Q(A,B)$. Moreover, if the walk is lazy then $\hat\Psi(A)=\hat\Q(A,A^c)=\Q(A,A^c)=\Psi(A)$, or if $\Psi(A)\leq\Delta_{min}/2$ (defined below) then again $\hat\Psi(A)=\Psi(A)$.

To motivate the form of our main result, we note that in their work on Blocking conductance Kannan, Lov\'asz and Montenegro \cite{KLM06.1} show that the square of conductance can often be replaced by a product of a measure of vertex boundary and a measure of edge expansion. Likewise, our general bound will involve a product of edge expansion $\hat\Psi_{min}$ with a measure of vertex boundary $\hat{A}_{max}$, rather than just the square of edge expansion which is found in most isoperimetric results.

\begin{theorem} \label{thm:general}
Given a finite Markov chain, let
$$
\begin{array}{lcllcl}
\vspace{1ex}\displaystyle
\hat\Psi_{min} &=& \displaystyle \min_{\pi(A)\leq 1/2}\hat\Psi(A) &
\hat{A}_{min} &=& \displaystyle \min\{\hat\Psi_{min},\,\Delta_{min}/2\} \\ \vspace{1ex}\displaystyle 
\Delta_{min} &=& \displaystyle \min_{\substack{A,B\subset V,\\ \pi(A)\neq\pi(B)}} |\pi(A)-\pi(B)| \qquad &
\hat{A}_{max} &=& \displaystyle \max\{\hat\Psi_{min},\,\Delta_{min}/2\} \\
\Q_{min} &=& \displaystyle \min_{A\subset V} \Q(A,A^c) &
\pi_* &=& \displaystyle \min_{x\in V}\pi(x) \,.
\end{array}
$$
Then,
\begin{eqnarray*}
\tau(\epsilon) &\leq& \frac{1}{-\log \left(1-2\frac{\hat{A}_{min}}{\Delta_{min}}\,(1-\cos(2\pi\,\hat{A}_{max}))\right)}\log \frac{1-\pi_*}{\epsilon} \\
 &\approx& \frac{1}{2\pi^2\hat\Psi_{min} \hat{A}_{max}} \log\frac{1-\pi_*}{\epsilon} \\ 
1-|\lambda_i| &\geq& 2\frac{\hat{A}_{min}}{\Delta_{min}}\,(1-\cos(2\pi\,\hat{A}_{max})) 
 \approx 2\pi^2\,\hat\Psi_{min}\,\hat{A}_{max} \\
\lambda &\geq& \frac{2\Q_{min}}{\pi_*}(1-\cos(\pi\,\pi_*))
 \approx \pi^2\,\pi_*\Q_{min}
\end{eqnarray*}
\end{theorem}

\begin{proof}[Proof of Collaries \ref{cor:Eulerian-simple} and \ref{cor:Eulerian-max} (see Section \ref{sec:distances} for the $L^{\infty}$-bounds)]
First to Corollary \ref{cor:Eulerian-simple}. Suppose that $\pi(A)\leq 1/2$, and $B\subset V$ satisfies $\Psi(A)=\Q(A,B)+\frac{\pi(A^c)-\pi(B)}{\pi(v)}\Q(A,v)$. If $N(A)\subseteq B^c$ then $\pi(N(A)\setminus v)\leq \pi(B^c\setminus v)=1-\pi(B\cup v)<\pi(A)$, contradicting the expansion condition. Hence, $N(A)\cap B\neq\emptyset$ and so $\exists x\in A,\ y\in B$ with $\P(x,y)>0$, and so $\Psi(A)\geq\Q(A,B)\geq\pi(x)\P(x,y)\geq 1/m$. Likewise, for some $B\subset V$, $\hat\Psi(A)\geq\hat\Q(A,B)\geq \min\{\Q(x,y),\pi(y)/2\}$ and so $\hat\Psi(A)\geq 1/m$ if $\pi(y)=deg(y)/m\geq 2/m$. If $deg(y)=1$ then $N(\{y\})$ has only a single vertex $v$, and so $\pi(N(\{y\})\setminus v)=0$ contradicting the expansion condition. It follows that $\pi(y)\geq\pi_*\geq 2/m$. Corollary \ref{cor:Eulerian-simple} then follows from Theorem \ref{thm:general} and the bound $\Delta_{min}\geq 1/m$. Corollary \ref{cor:Eulerian-max} follows similarly, but with $\hat\Psi(A)\geq 1/nd$ and $\Delta_{min}=1/n$. 
\end{proof}

Note that the max-degree walk is actually the same as the simple random walk when each vertex $x$ has $d-deg(x)$ self-loops added, and yet Corollary \ref{cor:Eulerian-max} is much better than that induced by Corollary \ref{cor:Eulerian-simple}. To understand this, recall that, in keeping with the intuition of Blocking Conductance, Theorem \ref{thm:general} will greatly improve on a bound involving edge-expansion alone (i.e. $\Psi(A)$ or $\hat\Psi(A)$) if $\Delta_{min}\gg \Psi_{min}$. In fact, the max-degree walk had $\Delta_{min}=1/n\gg 1/nd=\Psi_{min}$. 

The theorem gets us very close to answering the question of what is the worst of all random walks, as shown by the following examples.

\begin{example}
Consider the simple random walk on a cycle (with $m=2n$ edges). Since $m=2n$ and $d=2$ then Corollaries \ref{cor:Eulerian-simple} and \ref{cor:Eulerian-max} are the same. In Example \ref{ex:max-degree} we find that Corollary \ref{cor:Eulerian-max} is exact for the eigenvalue gap and essentially sharp for mixing times, and hence Corollary \ref{cor:Eulerian-simple} is equally good.
\end{example}

\begin{example} \label{ex:max-degree}
Consider a max-degree walk on a cycle with an odd number of vertices $n$, such that at each vertex there are $d-1$ edges pointing in the clockwise direction, and $1$ edge pointing in the counterclockwise direction.

This walk has an eigenvalue $\lambda_k=\frac{d-1}{d}e^{\pi i(n-1)/n} + \frac{1}{d}e^{-\pi i(n-1)/n}$ with eigenvector $f(x)=e^{\pi i x(n-1)/n}$  where the vertices are labeled clockwise as $x\in\{0,\,1,\,\ldots,\,n-1\}$. Then
$$
1-|\lambda_k|=1-\sqrt{1-\frac 4d(1-1/d)\sin^2\frac{\pi(n-1)}{n}}
 \approx \frac{2}{d}(1-1/d)\left(\frac{\pi}{n}\right)^2
 \approx \frac{2\pi^2}{n^2 d}\,.
$$
Corollary \ref{cor:Eulerian-max} gives a fairly similar bound of
$$
\min 1-|\lambda_k| \geq \frac{2}{d}(1-\cos(\pi/n)) \approx \frac{\pi^2}{n^2 d}\,.
$$
The upper and lower bounds are equal at $d=2$, and within a factor two of equality when $d>2$.

For spectral gap, note that $\frac{\P+\P^*}{2}$ is just the simple random walk on a cycle, and the largest eigenvalue of this is $\cos(2\pi/n)$. Consequently $\lambda=1-\cos(2\pi/n)\approx \frac{2\pi^2}{n^2}$. By Theorem \ref{thm:general} every walk with $\Q_{min}=1/n$ and $\pi_*=1/n$ satisfy$\lambda\geq 2(1-\cos(\pi/n))\approx \frac{\pi^2}{n^2}$, and so our drifting walk is within a factor two of having the worse spectral gap among all walks with $\Q_{min}=1/n$ and $\pi_*=1/n$. Although Corollary \ref{cor:Eulerian-max} is quite poor for this example, it is only off by a factor of four when considering instead the simple random walk on a cycle with $d-2$ self-loops (and $\lambda=1-\frac{2}{d}(1-\cos(2\pi/n))$.

Likewise, the upper and lower bounds on mixing time are quite similar:
$$
\tau(\epsilon) 
  \geq \frac{1}{-\log|\lambda_k|}\log \frac{1}{2\epsilon}
  \geq \frac{1}{-\log\sqrt{1-\frac 4d(1-1/d)\sin^2\frac{\pi(n-1)}{n}}}
        \log \frac{1}{2\epsilon} 
  \approx \frac{n^2 d}{2\pi^2}\log\frac{1}{2\epsilon}
$$
while the upper bound is
$$
\tau(\epsilon) 
  \leq \frac{1}{-\log(1-\frac 2d(1-\cos(\pi/n))}\log \frac{1-1/n}{\epsilon}
  \approx \frac{n^2 d}{\pi^2}\log\frac{1}{\epsilon}\,.
$$
The bounds are nearly equivalent at $d=2$, and within a factor two of equality when $d>2$. When $n=3$ and $d=2$ then the lower bound can be sharpened slightly to be exactly equal to the upper bound.
\end{example}

\begin{example}
Consider a general Markov chain. Note that if vertex $v\in V$ has $\pi(v)=\pi_*$ then $\Psi_{\min}\leq\Psi(\{v\})\leq\pi_*(1-\pi_*)$. The clockwise Markov chain given above had $\Psi_{min}=1/nd=\pi_*/d$. If instead the walk has transitions $\P(x,x+1)=\alpha\in[1/2,1]$ and $\P(x,x-1)=1-\alpha$ then $\Psi_{min}=(1-\alpha)\pi_*$. If $\Psi_{min}\leq\pi_*/2$ then when $\alpha=1-\Psi_{min}/\pi_*\geq 1/2$ the upper and lower bounds in Theorem \ref{thm:main} are again within a factor $2$ or $4$ from the correct values.
\end{example}

\begin{proof}[Proof of Theorem \ref{thm:general}]
As suggested in the preliminaries, we will study the $f$-congestion $\C_{\sin(\pi\,a)}$, via Lemma \ref{lem:worst-case}. This will be done in two steps. First, we show a result appropriate for max-degree random walks. Then we consider a case relevant to the simple random walk.

Fix some set $A\subset V$. 

First consider the case that $\Psi(A)<\Delta_{min}/2$. 

Notice that if $\pi(A_u)>\pi(A)$ then $\pi(A_u)\geq\pi(A)+\Delta_{min}$, while if $\pi(A_u)<\pi(A)$ then $\pi(A_u)\leq\pi(A)-\Delta_{min}$. In Figure \ref{fig:uniform} let the solid line sketch the curve $\pi(A_u)$, and note that the dashed line $m(u)$ encloses the same area $\Psi(A)$ but decreases the integral, so $\int_0^t m(u)\,du \leq \int_0^t \pi(A_u)\,du$.

\begin{figure}[ht]
\begin{center}
\includegraphics[height=1.5in]{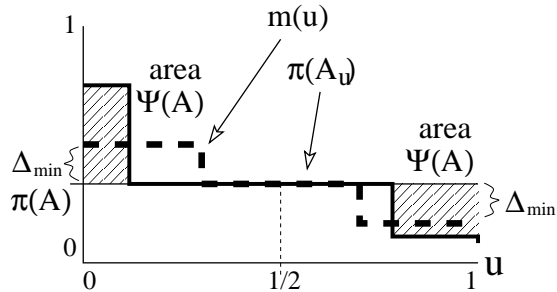}
\caption{If $\Psi(A)<\frac{\Delta_{min}}{2}$ then $\int_0^t m(u)\,du \leq \int_0^t \pi(A_u)\,du$ and $\int_0^1 m(u)\,du=\int_0^1\pi(A_u)\,du$.} \label{fig:uniform}
\end{center}
\end{figure}

The curve $m(u)$ is given by
$$
m(u) = 
\begin{cases}
\pi(A)+\Delta_{min} & \textrm{if $u\leq \Psi(A)/\Delta_{min}$} \\
\pi(A)              & \textrm{if $u \in(\Psi(A)/\Delta_{min},1-\Psi(A)/\Delta_{min})$} \\
\pi(A)-\Delta_{min} & \textrm{if $u\geq 1-\Psi(A)/\Delta_{min}$}\,.
\end{cases}
$$

By Lemma \ref{lem:worst-case}, since $\int_0^1 m(u)\,du=\int_0^1\pi(A_u)\,du=\pi(A)$, it follows that if $\pi(A)=x$ then
\begin{eqnarray*}
\C_{\sin(\pi a)}(A) &\leq& \frac{\int_0^1 \sin(\pi m(u))\,du}{\sin(\pi x)} \\
 &=& 1-2\frac{\Psi(A)}{\Delta_{min}}\left(1-\cos(\pi\Delta_{min})\right) \\
 &\leq& 1-2\frac{\hat\Psi(A)}{\Delta_{min}}\left(1-\cos(\pi\Delta_{min})\right) \,.
\end{eqnarray*}

Now, consider the case that $\Psi(A)\geq\Delta_{min}/2$.

In Figure \ref{fig:general} let the solid line in the left diagram sketch the curve $\pi(A_u)$. If this has portions above $y=\pi(A)$ but beyond $u=1/2$ then truncate these off, with the portion below $y=\pi(A)$ raised slightly to keep the enclosed area constant, making the left figure into the right side one. Note that the dashed line $m(u)$ encloses the same area $\hat\Psi(A)$ but decreases the integral, so $\int_0^t m(u)\,du \leq \int_0^t \pi(A_u)\,du$ and $\int_0^1 m(u)\,du=\int_0^1\pi(A_u)\,du=\pi(A)$.

\begin{figure}[ht]
\begin{center}
\includegraphics[height=1.5in]{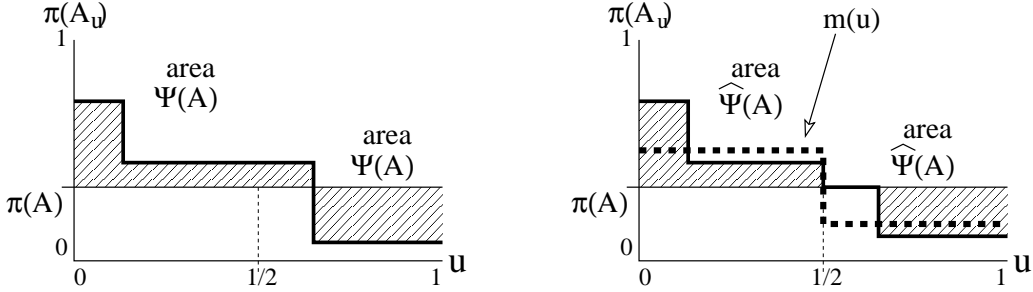}
\caption{Curve $m(u)$ satisfies $\int_0^t m(u)\,du\leq\int_0^t\pi(A_u)\,du$ and $\int_0^1 m(u)\,du=\int_0^1 \pi(A_u)\,du$.} \label{fig:general}
\end{center}
\end{figure}

The curve $m(u)$ is such that, if $\wp\in[0,1/2]$ denotes the value of $u$ where $m(u)$ crosses the line $y=\pi(A)$, that is $m(u)>\pi(A)$ if $u<\wp$ and $m(u)<\pi(A)$ if $u>\wp$, then
$$
m(u) = 
\begin{cases}
\pi(A)+\frac{\hat\Psi(A)}{\wp} & \textrm{if $u\leq \wp$} \\
\pi(A)-\frac{\hat\Psi(A)}{1-\wp} & \textrm{if $u>\wp$}
\end{cases}
$$

By Lemma \ref{lem:worst-case} it follows that, if $g(u)=\pi(A_u)$ and $x=\pi(A)$ then
\begin{eqnarray*}
\int_0^1 \sin(\pi\,g(u))\,du &\leq& \int_0^1 \sin(\pi\,m(u))\,du \\
  &=& \wp\,\sin\left(\pi\left(x + \frac{\hat\Psi(A)}{\wp}\right)\right) 
  + (1-\wp)\,\sin\left(\pi\left(x - \frac{\hat\Psi(A)}{1-\wp}\right)\right) \\
  &\leq& \sin(\pi x)\,\cos(2\pi\hat\Psi(A))\,,
\end{eqnarray*}
where the final inequality is from Lemma \ref{lem:ineq} in the Appendix.
Hence, $\C_{\sin(\pi a)}(A)\leq\cos(2\pi\hat\Psi(A))$.

Combine these two cases, maximize over sets $A\subset V$, and apply Theorem \ref{thm:main} to obtain the mixing time and eigenvalue bounds.

For the spectral gap, note that
$$
\lambda = 2\min 1-\lambda_i\left(\frac{\I}{2}+\frac{\P+\P^*}{4}\right) 
  = 2\min 1-\left|\lambda_i\left(\frac{\I}{2}+\frac{\P+\P^*}{4}\right)\right|\,.
$$
This shows it suffices to study eigenvalues of $\P'=\frac{\I}{2}+\frac{\P+\P^*}{4}$. However, $\P'$ is a lazy walk, and so $\hat\Psi(A)=\Psi(A)=\Q_{\P'}(A,A^c)$. This is in turn half the ergodic flow $\Q_{\frac{\P+\P^*}{2}}(A,A^c)$, and so $\hat\Psi(A)=\frac 12\Q_{\frac{\P+\P^*}{2}}(A,A^c)=\frac 12\Q_{\P}(A,A^c)$ (since $\Q_{\P}(A,A^c)=\Q_{\P^*}(A,A^c)=\Q_{\frac{\P+\P^*}{2}}(A,A^c)$). In short,
$$
\hat\Psi_{min}\left(\frac{\I}{2}+\frac{\P+\P^*}{4}\right) 
  = \frac 12\min_{A\subset V}\Q_{\P}(A,A^c) = \frac 12\Q_{min}(\P)\,.
$$
Before applying the eigenvalue bounds proven earlier, note for a lazy walk that $\Psi(A)=\Q(A,A^c)$, with $\frac{\Q(A,x)}{\pi(x)}< \frac 12$ only if $x\in A^c$, and $\frac{\Q(A,x)}{\pi(x)}>\frac 12$ only if $x\in A$. It follows that if $\pi(A_u)>\pi(A)$ then $A\subsetneq A_u$ and so $\pi(A_u)\geq\pi(A)+\pi_*$. Likewise, if $\pi(A_u)<\pi(A)$ then $A\supsetneq A_u$ and so $\pi(A_u)\leq\pi(A)-\pi_*$. Hence, when studying a lazy walk (such as $\frac{\I}{2}+\frac{\P+\P^*}{4}$), $\Delta_{min}$ may be replaced by $\pi_*$ in our earlier analysis. But $\pi_*\left(\frac{\I}{2}+\frac{\P+\P^*}{4}\right)=\pi_*(\P)$, and so the spectral bound follows from the earlier eigenvalue bounds.
\end{proof}


\section{Other distances} \label{sec:distances}

Total variation distance measures only distance from stationary at an average vertex. The much stronger $L^{\infty}$ distance measures distance from stationary at the worst vertex. In this section we show mixing bounds on $L^{\infty}$, $L^2$ and relative entropy distances which are again within a small constant factor of those for the walk on a cycle with clockwise drift.

Given distributions $\sigma$ and $\pi$, the relative entropy distance ${\rm D}(\sigma\|\pi)$, the $L^2$ distance $\|\sigma/\pi-1\|_{2,\pi}$ and the relative pointwise (or $L^{\infty}$) distance $\|\sigma/\pi-1\|_{\infty,\pi}$ are defined by
\begin{eqnarray*}
{\rm D}(\sigma\|\pi) &=& \sum_{x\in V} \pi(x)\left(\frac{\sigma(x)}{\pi(x)}\right)\log \frac{\sigma(x)}{\pi(x)} \\
\|\sigma/\pi-1\|_{2,\pi} &=& \sqrt{\sum_{x\in V} \pi(x)\left(\frac{\sigma(x)}{\pi(x)}-1\right)^2} \\
\|\sigma/\pi-1\|_{\infty} &=& \max_{x\in V}\left|\frac{\sigma(x)}{\pi(x)}-1\right|
\end{eqnarray*}
The worst case number of steps required for a walk to reach distance $\epsilon$ is given by the mixing times $\tau_D(\epsilon)$, $\tau_2(\epsilon)$ and $\tau_{\infty}(\epsilon)$ respectively. These are related by $\tau(\epsilon)\leq \tau_D(2\epsilon^2)$, $\tau(\epsilon)\leq\tau_2(2\epsilon)$, $\tau_D(\epsilon)\leq \tau_2(\sqrt{\epsilon})$, $\tau_2(\epsilon)\leq\tau_{\infty}(\epsilon) \leq \tau_2\left(\epsilon\sqrt{\frac{1-\pi_*}{\pi_*}}\right)$, and $\tau_{\infty}(\epsilon)\leq\tau_{2,\P}(\sqrt{\epsilon})+\tau_{2,\P^*}(\sqrt{\epsilon})$ where $\tau_{2,K}(\epsilon)$ denotes $L^2$ mixing for the Markov chain $K$ (e.g. see Appendix of \cite{MT06.1} for proofs and/or references for these).

Evolving set bounds on relative entropy and $L^2$ mixing (see \cite{Mon06.1} or \cite{MT06.1}) can be combined with the technique used to prove Theorem \ref{thm:main} (see \cite{Mon06.4}), to show
\begin{eqnarray*}
\tau_D(\epsilon) &\leq& \frac{1}{1-\C_{\sin(\pi a)}}\left(\log \log \frac{1}{\pi_*}+\log\frac{1}{\epsilon}\right) \\
\tau_2(\epsilon) &\leq& \frac{1}{1-\C_{\sin(\pi a)}}\left(\frac 12\log \frac{1-\pi_*}{\pi_*} + \log \frac{1}{\epsilon}\right) \\
\tau_{\infty}(\epsilon) &\leq& \frac{1}{1-\C_{\sin(\pi a)}}\left(\log \frac{1-\pi_*}{\pi_*} + \log \frac{1}{\epsilon}\right) 
\end{eqnarray*}
The $L^{\infty}$ bound followed from the $L^2$ bound and the relation $\tau_{\infty}(\epsilon)\leq\tau_2\left(\epsilon\sqrt{\frac{1-\pi_*}{\pi_*}}\right)$.

When $\epsilon\to 0$ these are asymptotically the same as the $\tau(\epsilon)$ bound of Theorem \ref{thm:general}. However, when $\epsilon$ is large we can further improve these via Evolving set bounds shown in Montenegro \cite{Mon06.1} (see also \cite{MT06.1}):

\begin{theorem} \label{thm:lots-o-evosets}
For a finite Markov chain, if $\C_f(r)=\max_{\pi(A)\leq r}\C_f(A)$ $\forall r\in[0,1]$ then 
\begin{eqnarray*}
\tau(\epsilon)
  &\leq& \left\lceil \int_{\pi_*}^{1-\epsilon} \frac{dr}{(1-r)(1-\C_{a(1-a)}(r))} \right\rceil \\
\tau_D(\epsilon)
  &\leq& \left\lceil \int_{\pi_*}^{e^{-\epsilon}} \frac{dr}{r\log(1/r)(1-\C_{a\log(1/a)}(r))} \right\rceil \\
\tau_2(\epsilon)
  &\leq& \left\lceil \int_{\pi_*}^{1/(1+\epsilon^2)} \frac{dr}{2r(1-r)(1-\C_{\sqrt{a(1-a)}}(r))} \right\rceil
\end{eqnarray*}
where the total variation bound requires $r\left(1-\C_{a(1-a)}(1-r)\right)$ be convex, relative entropy requires $r\left(1-\C_{a\log(1/a)}(e^{-r})\right)$ to be convex, and the $L^2$ bound requires $r\left(1-\C_{\sqrt{a(1-a)}}\left(\frac{1}{1+r^2}\right)\right)$ to be convex.
\end{theorem}

It remains to bound $\C_f(r)$ for each choice of $f(a)$ given above, then integrate in Theorem \ref{thm:lots-o-evosets}. First, the bounds on $f$-congestion for the quantities used in the various distances of interest to us.

\begin{lemma} \label{lem:lots-o-mixing}
\begin{eqnarray*}
1-\C_{a(1-a)}(r)        
  &\geq& \frac{2\Delta_{min}\Psi_{min}}{r(1-r)}\,\delta_{r\leq 1/2} + 8\Delta_{min}\Psi_{min} \\
1-\C_{a\log(1/a)}(r)    
  &\geq& \frac{\Delta_{min}\Psi_{min}}{2r^2\log \frac 1r}\,\delta_{r\leq e^{-1/2}} + e\Delta_{min}\Psi_{min} \\
1-\C_{\sqrt{a(1-a)}}(r) 
  &\geq& \frac{\Delta_{min}\Psi_{min}}{4r^2(1-r)^2}\,\delta_{r\leq 1/2} + 4\Delta_{min}\Psi_{min}\,\delta_{r>1/2}
\end{eqnarray*}
\end{lemma}

\begin{proof}
Montenegro \cite{Mon06.1} uses Lemma \ref{lem:worst-case}, and optimization similar to the use of Lemma \ref{lem:ineq} in the proof of Theorem \ref{thm:main}, to lower bound the $f$-congestion quantities in terms of $\Psi(A)$:
$$
\begin{array}{lclcl}
\vspace{1ex}\displaystyle 1-\C_{a(1-a)}(r)        
  &\geq& \displaystyle \min_{\pi(A)\leq r} \frac{4\Psi(A)^2}{\pi(A)\pi(A^c)} 
  &\geq& \displaystyle \frac{4\Psi_{min}^2}{r(1-r)}\,\delta_{r\leq 1/2}+16\Psi_{min}^2\,\delta_{r>1/2} \\
\vspace{1ex}\displaystyle 1-\C_{a\log(1/a)}(r)  
  &\geq& \displaystyle \min_{\pi(A)\leq r} \frac{2\Psi(A)^2}{\pi(A)^2\log \frac{1}{\pi(A)}} 
  &\geq& \displaystyle \frac{2\Psi_{min}^2}{r^2\log\frac 1r}\,\delta_{r\leq e^{-1/2}}+4e\Psi_{min}^2\,\delta_{r>e^{-1/2}} \\
\displaystyle 1-\C_{\sqrt{a(1-a)}}(r) 
  &\geq& \displaystyle \min_{\pi(A)\leq r} \frac{\Psi(A)^2}{2\pi(A)^2\pi(A^c)^2}
  &\geq& \displaystyle \frac{\Psi_{min}^2}{2r^2(1-r)^2}\,\delta_{r\leq 1/2} + 8\Psi_{min}^2\,\delta_{r>1/2}
\end{array}
$$

Improvements are again possible for max-degree type walks. If $\Psi_{min}\leq\Delta_{min}/2$ then the worst case $m(u)$ was already determined in the proof of Theorem \ref{thm:general} (recall that Lemma \ref{lem:worst-case} says $m(u)$ is the worst for every concave function $f:[0,1]\to\R_+$). Notice that the upper bound on $\C_{\sin(\pi a)}(A)$ could have been obtained by taking 
$$
1-\C_{f(a)}(A)\geq 2\frac{\Psi(A)}{\Delta_{min}}\left(1-\left.\C_{f(a)}(\pi(A))\right|_{\Psi(A)=\Delta_{min}/2}\right)\,.
$$
In short, the case with $\Psi_{min}\leq\Delta_{min}/2$ can be reduced to the case of $\Psi_{min}\geq\Delta_{min}/2$, without losing any accuracy. However, bounds in terms of $\Psi_{min}$ were found at the beginning of this proof, so it suffices to replace $\Psi(A)$ by $\Delta_{min}/2$ in these bounds, and then multiply by $\frac{2\Psi(A)}{\Delta_{min}}$.

Combining the two cases $\Psi(A)>\Delta_{min}/2$ and $\Psi(A)\leq\Delta_{min}/2$ gives the lemma.
\end{proof}

Applying Theorems \ref{thm:lots-o-evosets} and \ref{thm:main} to the $f$-congestion bounds of Lemma \ref{lem:lots-o-mixing} leads to the following relations:

\begin{corollary} \label{cor:lots-o-mixing}
For a finite Markov chain, if $\Psi_{min}=\displaystyle\min_{A\subset V}\Psi(A)$ and $A_{max}=\max\{\Psi_{min},\frac{\Delta_{min}}{2}\}$ then
\begin{eqnarray*}
1-|\lambda_i| &\geq& 16\Psi_{min}A_{max} \\
\tau(\epsilon) &\leq& \left\lceil \frac{\frac{1-4\pi_*^2}{2}+\log\frac{1}{2\epsilon}}{16\Psi_{min}A_{max}}\,\delta_{\epsilon\leq 1/2}
  + \frac{(1-\epsilon)^2-\pi_*^2}{8\Psi_{min}A_{max}}\,\delta_{\epsilon>1/2} \right\rceil \\
\tau_D(\epsilon) &\leq& \left\lceil \frac{1-e\pi_*^2+\log\frac{1}{2\epsilon}}{4e\Psi_{min}A_{max}}\,\delta_{\epsilon\leq 1/2} 
  + \frac{e^{-2\epsilon}-\pi_*^2}{4\Psi_{min}A_{max}}\,\delta_{\epsilon> 1/2} \right\rceil \\
\tau_2(\epsilon) &\leq& \left\lceil \frac{\frac 23+\log\frac{1}{\epsilon}}{8\Psi_{min}A_{max}}\,\delta_{\epsilon\leq 1}
  + \frac{1+3\epsilon^2}{6\Psi_{min}A_{max}(1+\epsilon^2)^3}\,\delta_{\epsilon>1} \right\rceil 
\end{eqnarray*}
\end{corollary}

Note that $\Psi_{min}\geq\hat\Psi_{min}$, and so even the total variation bound can occasionally improve on Theorem \ref{thm:general}. For instance, the simple random walk on an expanding Eulerian graph has $\Psi_{min}\geq 1/m$, and so the bound here is better than that of Corollary \ref{cor:Eulerian-simple} when $m$ is large and $\epsilon>0.361$. This explains why our upper bound on mixing time was not quite sharp before.

\begin{proof}[Proof of $L^{\infty}$ cases in Corollaries \ref{cor:Eulerian-simple} and \ref{cor:Eulerian-max}]
We use the relation $\tau_{\infty}(\epsilon)\leq\tau_{2,\P}(\sqrt{\epsilon})+\tau_{2,\P^*}(\sqrt{\epsilon})$.

First, Corollary \ref{cor:Eulerian-max}. If $A,B\subset V$ and $\pi(B)=\pi(A^c)$ then
\begin{equation} \label{eqn:reversing}
\Q_{\P}(A,B)=\pi(B)-\Q_{\P}(A^c,B)=\pi(B)-\pi(A^c)+\Q_{\P}(A^c,B^c)=\Q_{\P^*}(B^c,A^c)\,.
\end{equation}
Since $\pi$ is uniform for the max-degree walk then $\Psi(A)=\min_{\pi(B)=\pi(A^c)}\Q(A,B)$, and so by \eqref{eqn:reversing} $\Psi_{min}$ is the same for $\P$ and $\P^*$. Then $\tau_{\infty}(\epsilon)\leq 2\tau_2(\sqrt{\epsilon})$ if the $L^2$ mixing bound from Corollary \ref{cor:lots-o-mixing} is used, that is, if we substitute the conditions $\Delta_{min}\geq 1/n$, $\Psi_{min}\geq 1/nd$ and $A_{max}\geq 1/2n$.

Now to Corollary \ref{cor:Eulerian-simple}. If $A,B\subset V$, $v\in V$ and $\pi(B)\leq\pi(A^c)<\pi(B\cup v)$ then, arguing as in \eqref{eqn:reversing},
\begin{equation} \label{eqn:reversing2}
\Q_{\P^*}(A,B) + \frac{\pi(A^c)-\pi(B)}{\pi(v)}\,\Q_{\P^*}(A,v) 
  = \Q_{\P}(B^c\setminus v,A^c) + \left(1-\frac{\pi(A^c)-\pi(B)}{\pi(v)}\right)\,\Q_{\P}(v,A^c)\,.
\end{equation}
Arguing as in the proof of Corollaries \ref{cor:Eulerian-simple} and \ref{cor:Eulerian-max} after Theorem \ref{thm:general}, if $\pi(C)\leq 1/2$ and $\pi(D)\leq\pi(C^c)<\pi(D\cup v)$ then $\Q(C,D)\geq 1/m$. Note that for some $A,B,v$ with $\pi(A)\leq 1/2$, equation \eqref{eqn:reversing2} is exactly $\Psi_{min}$ for the walk $\P^*$, and so if $C=B^c\setminus v$ and $D=A^c$ then it follows that $\Psi_{min}\geq \Q_{\P}(B^c\setminus v,A^c)\geq 1/m$. Hence $\Delta_{min}\geq 1/m$, $\Psi_{min}\geq 1/m$ and $A_{max}\geq 1/2m$ for both $\P$ and $\P^*$. The $L^{\infty}$ case in Corollary \ref{cor:Eulerian-simple} follows by the resulting bounds on $\tau_{2,\P}(\sqrt{\epsilon})$ and $\tau_{2,\P^*}(\sqrt{\epsilon})$ in Corollary \ref{cor:lots-o-mixing}. 
\end{proof}


\bibliographystyle{plain}
\bibliography{../references}


\section*{Appendix}

In the Appendix we look at two results needed in this paper. First, the proof of Theorem \ref{thm:main}, and then the proof of an inequality used in showing Theorem \ref{thm:general}.

\begin{proof}[Proof of Theorem \ref{thm:main}]
Given $x\in V$, the Evolving set process is defined recursively by setting $S_0=\{x\}$, and then to determine $S_{t+1}$ choose $u\in[0,1]$ uniformly at random, and set $S_{t+1}=(S_t)_u$. Let $\E_t$ denote the expectation after $t$ steps of the Evolving set process. Also, we use the notation $S^{\#}$ to denote $S$ if $\pi(S)\leq 1/2$ and $S^c=V\setminus S$ if $\pi(A)>1/2$. 

Starting with an inequality of \cite{Mon06.1} (see also \cite{MT06.1}), we have
\begin{eqnarray*}
\|\P^t(x,\cdot)-\pi\|_{TV} &\leq& \frac{1}{\pi(x)}\,\E_t \pi(S_t)(1-\pi(S_t)) \\
   &\leq& \left(\max_{\substack{A\subset V,\\ \pi(A)\leq 1/2}} \frac{\pi(A)\pi(A^c)}{f(\pi(A))}\right) \frac{1}{\pi(x)}\,\E_t\,f(\pi(S_t^{\#})) \\
   &\leq& \left(\max_{\substack{A\subset V,\\ \pi(A)\leq 1/2}} \frac{\pi(A)\pi(A^c)}{f(\pi(A))}\right)\frac{1}{\pi(x)}\,\E_{t-1}\,f(\pi(S_{t-1}^{\#}))\,\C_f(S_{t-1}^{\#}) \\
   &\leq& \left(\max_{\substack{A\subset V,\\ \pi(A)\leq 1/2}} \frac{\pi(A)\pi(A^c)}{f(\pi(A))}\right)\frac{f(\pi(\{x\}^{\#}))}{\pi(x)}\,\C_f^t\,.
\end{eqnarray*}
The final inequality followed from $\C_f(S_{t-1}^{\#})\leq\C_f$, and then induction.

The lower bound on total variation distance can be found in \cite{Mon06.2} (see also \cite{MT06.1}).

The bound on eigenvalues follows by combining the upper and lower bounds:
$$
|\lambda_i| \leq \sqrt[t]{2\C_f^t\max_{\substack{x\in V,\\\pi(A)\leq 1/2}} \frac{\pi(A)\pi(A^c)}{f(\pi(A))}\frac{f(\pi(\{x\}^{\#}))}{\pi(x)}}\xrightarrow{t\to\infty} \C_f\,.
$$
\end{proof}

We have left for the Appendix the proof of an inequality key to our main theorem.

\begin{lemma} \label{lem:ineq}
If $0\leq a \leq \frac 12$, $\frac{c}{1-a}\leq b\leq 1/2$ and $0\leq c\leq b(1-b)$ then
\begin{eqnarray*}
h(a,b,c) &=& b\sin\left(\pi\left(a+\frac{c}{b}\right)\right) + (1-b)\sin\left(\pi\left(a-\frac{c}{1-b}\right)\right) \\
 &\leq& \sin(\pi a)\cos(2\pi c)
\end{eqnarray*}
%
\end{lemma} 

\begin{proof}
First, use the expansion $\sin(x+y)=\sin(x)\cos(y)+\sin(y)\cos(x)$ to re-arrange the terms a bit:
\begin{eqnarray}
\lefteqn{h(a,b,c)} \label{eqn:expansion} \\
&=& \sin(\pi a)\,\left[\left(b\cos \frac{\pi c}{b} + (1-b)\cos\frac{\pi c}{1-b}\right) + \cot(\pi a)\,\left(b\sin \frac{\pi c}{b} - (1-b)\sin\frac{\pi c}{1-b}\right)\right]
\nonumber
\end{eqnarray}

Consider the second term. Suppose $b\in[c,1]$. Then
\begin{eqnarray*}
\frac{d}{db} b\sin\frac{\pi c}{b} &=& \sin\frac{\pi c}{b} - \frac{\pi c}{b}\cos \frac{\pi c}{b} \\
 &=& \cos\frac{\pi c}{b}\left(\tan\frac{\pi c}{b}-\frac{\pi c}{b}\right) \\
 &\geq& 0
\end{eqnarray*}
The inequality is because $(\cos x)(\tan x-x)\geq 0$ when $x=\frac{\pi c}{b}\in[0,\pi]$. 

It follows that $b\sin\frac{\pi c}{b}$ is increasing, and in particular if $b\leq 1/2$ then
$$
b\sin \frac{\pi c}{b} - (1-b)\sin\frac{\pi c}{1-b} \leq 0\,.
$$
Consequently, if $b$ and $c$ are fixed then the ratio $\frac{h(a,b,c)}{\sin(\pi a)}$ is maximized when $a\in[0,1/2]$ is maximized. Subject to the conditions on $a$, $b$ and $c$ in the lemma this maximum is at $a=\min\{1/2,1-c/b\}$.

It has just been shown that if $b\geq 2c$ then $\frac{h(a,b,c)}{\sin(\pi a)} \leq \frac{h(1/2,b,c)}{\sin(\pi/2)}=h(1/2,b,c)$, otherwise $\frac{h(a,b,c)}{\sin(\pi a)} \leq \frac{h(1-c/b,b,c)}{\sin(\pi(1-c/b))}$. The latter case can be simplified further by the relation $\frac{h(1-c/b,b,c)}{\sin(\pi(1-c/b))}\leq h(1/2,2c,c)$ when $b< 2c$. This requires showing that
$$
(1-b)\sin \frac{\pi c}{b(1-b)} \leq (1-2c)\sin \frac{\pi c}{b}\,\cos \frac{\pi c}{1-2c}
$$
where $\frac{b}{2} \leq c \leq b(1-b)$. The substitutions $x=\frac{\pi c}{b(1-b)}$ and $y=\frac{\pi c}{b}$  (ie. $b=1-y/x$ and $c=\frac{y(1-y/x)}{\pi}$) reduce this to Lemma \ref{lem:sinc}, and so the inequality holds.

Combining results, we now know that $\frac{h(a,b,c)}{\sin(\pi a)}\leq h(1/2,\min\{2c,b\},c)$. The second term in Equation \eqref{eqn:expansion} is zero when $a=1/2$, and so to study $h(1/2,\max\{2c,b\},c)$ it suffices to consider the first term in Equation \eqref{eqn:expansion}. Note that
\begin{eqnarray*}
\frac{d}{db} b\cos\frac{\pi c}{b} 
  &=& \cos\frac{\pi c}{b}   + \frac{\pi c}{b}\sin \frac{\pi c}{b} \\
\frac{d^2}{db^2} b\cos\frac{\pi c}{b}
  &=& -\frac{(\pi c)^2}{b^3}    \cos \frac{\pi c}{b}
\end{eqnarray*}
The second derivative is negative when $b\geq 2c$, and so $g(b,c)=b\cos\frac{\pi c}{b}$ is concave in $b\geq 2c$. It follows that if $1/2\geq b'>\max\{2c,b\}$ then $h(1/2,\max\{2c,b\},c)=g(\max\{2c,b\},c)+g(1-\max\{2c,b\},c)\leq g(b',c)+g(1-b',c)=h(1/2,b',c)$, and in particular,
\begin{equation}\label{eqn:half-time}
h(1/2,\max\{2c,b\},c)\leq h(1/2,b',c) \leq h(1/2,1/2,c)\,.
\end{equation}
The result then follows. 
\end{proof}

The following Lemma was required in the preceding proof.

\begin{lemma} \label{lem:sinc}
If $\pi/2\leq y<x\leq\pi$ and ${\rm sinc}(z):=\frac{\sin z}{z}$ then
$$
{\rm sinc}(y)\,\left(1-\frac{2}{\pi}y\left(1-\frac yx\right)\right)\cos\frac{y\left(1-\frac yx\right)}{1-\frac{2}{\pi}y\left(1-\frac yx\right)}
  \geq {\rm sinc}(x)\,.
$$
\end{lemma}

The function ${\rm sinc}(y)$ is decreasing when $y\in[\pi/2,\pi]$, and so the lemma says that the term after ${\rm sinc}(y)$ is a measure of how much the function drops between $y$ and $x$. A slightly weaker result that is perhaps a bit easier to grasp is
$$
\forall x,y\in[\pi/2,\pi],\,x>y:\,{\rm sinc}(x) \leq {\rm sinc}(y)\,\cos\left(2y(1-y/x)\right)\,.
$$

\begin{proof}
Rewrite the problem as
$$
f(x,y) := \left(1-\frac{2}{\pi}y\left(1-\frac yx\right)\right)\cos\frac{y\left(1-\frac yx\right)}{1-\frac{2}{\pi}y\left(1-\frac yx\right)}
  - \frac{\sin x}{x}\frac{y}{\sin y} \geq 0\,.
$$

Observe that $\forall x\in[\pi/2,\pi]:\,f(x,x)=f(x,\pi/2)=0$, and so the lemma holds at the two extreme values for $y$, that is $y=\pi/2$ and $y=x$. Moreover, given fixed $x\in[\pi/2,\pi]$, the first partial with respect to $y$ at $y=\pi/2$ is positive:
\begin{eqnarray*}
\left.\frac{\partial f}{\partial y}\right|_{y=\pi/2} 
  &=& \frac{2x(\pi-x)\cos(\pi-x)-(1-2(\pi-x))\sin(\pi-x)}{\pi x} \\
  &\geq& \min_{x\in[\pi/2,\pi]} \left\{0,\frac{2x(\pi-x)(1-\frac{2}{\pi}(\pi-x))-(1-2(\pi-x))(\pi-x)}{\pi x}\right\} \\
  &=& \min_{x\in[\pi/2,\pi]} \left\{0,\,\frac{\pi-x}{\pi x}\left(2\pi-1-\frac{4x(\pi-x)}{\pi}\right)\right\} = 0\,.
\end{eqnarray*}
The first inequality is because the expression is trivially positive if $\pi-x\in[1/2,\pi/2]$, whereas if $a:=\pi-x\in[0,1/2]$ then use the relations $\cos a \geq 1-\frac{2}{\pi}a$ and $\sin a\leq a$. It follows that the inequality $f(x,y)\geq 0$ also holds near $y=\pi/2$.

Now, consider the third partial derivative with respect to $y$:
\begin{eqnarray*}
\frac{\partial^3 f}{\partial y^3}
  &=& -\frac{\pi^5 x^2 (2y-x)^3}{(\pi x-2xy+2y^2)^5}\,\sin\frac{y\left(1-\frac yx\right)}{1-\frac{2}{\pi}y\left(1-\frac yx\right)} \\
  &&  -\frac{6\pi^3x(2y-x)(x(\pi-x)+2y(x-y))}{(\pi x-2xy+2y^2)^4}\,\cos\frac{y\left(1-\frac yx\right)}{1-\frac{2}{\pi}y\left(1-\frac yx\right)} \\
  &&  +\frac{6y}{x}\frac{\sin x}{\sin y}\cot^3 y -\frac{6}{x}\frac{\sin x}{\sin y}\cot^2 y +\frac{5y}{x}\frac{\sin x}{\sin y}\cot y - \frac{3}{x}\frac{\sin x}{\sin y}
\\
  &\leq& 0
\end{eqnarray*}
The third derivative was negative because every term in it is negative when $\pi\geq x\geq y\geq\pi/2$ (note that $\frac{y(1-y/x)}{1-\frac{2}{\pi}y(1-y/x)}\in[0,\pi/2]$ and $\pi x-2xy+2y^2>0$). 

From the third partial we know that the second partial is decreasing, and so for each $x\in[\pi/2,\pi]$ there are three possible cases: strictly convex in $y$, convex then concave in $y$, or strictly concave in $y$. The function cannot be strictly convex because it is zero at $y=\pi/2$, then increasing, and a convex function could not then be zero again at $y=x$. If it is convex then concave, then the convex portion is strictly increasing because $\frac{\partial f}{\partial y}(x,\pi/2)>0$, while the concave portion starts at a positive value and ends at $f(x,x)=0$. Hence the minimum of the convex portion is at $y=\pi/2$, and the minimum of the concave portion is at $y=x$. Finally, if it is strictly concave then the minimum is at an endpoint, so either $f(x,\pi/2)$ or $f(x,x)$. 

It follows that for each value of $x\in[\pi/2,\pi]$ the minimum is either $f(x,\pi/2)=0$ or $f(x,x)=0$. Hence the function is non-negative.
\end{proof}

\end{document}